\magnification=1200
\def\uom{\Omega}
\def\lam{\lambda}
\def\IR{\hbox{\rm I\kern-.2em\hbox{\rm R}}}
\let\mathcal=\cal
\def\sqr#1#2{{\vcenter{\vbox{\hrule height.#2pt
     \hbox{\vrule width.#2pt height #1pt \kern#1pt
     \vrule width.#2pt}
     \hrule height.#2pt}}}}
\def\square{\mathchoice\sqr56\sqr56\sqr{2.1}3\sqr{1.5}3}

\vskip1.4cm \noindent


\centerline 
{\bf On the fundamental eigenvalue ratio of the $p$--Laplacian}
\bigskip


{\baselineskip=.5ex
\vfootnote{}{\sevenrm
\noindent\copyright 2004 by the authors.
Reproduction of this article, in its entirety, by any means is permitted
for non--commercial purposes.}}

\centerline
{{\bf Jacqueline Fleckinger}\footnote{$^\heartsuit$}{{\tt jfleck@univ-tlse1.fr},
CEREMATH \& UMR MIP, Universit\'e Toulouse-1,21 all\'ees de Brienne, 
31000 Toulouse, France}, 
{\bf Evans M. Harrell II}
\footnote {$^\clubsuit$}{{\tt harrell@math.gatech.edu}, School of Mathematics,
Georgia Tech, Atlanta, GA\break 30332-0160, USA.  
This work is supported by NSF grant DMS-0204059},
{ \bf Fran\c cois de Th\'elin}
\footnote{$^\spadesuit$}{{\tt dethelin@mip.ups-tlse.fr},
UMR MIP, Universit\'e Paul Sabatier, 31062 Toulouse, France}}

\vskip .7 true in

\centerline{\bf Abstract}

\smallskip
It is shown that the fundamental eigenvalue ratio
$\lam_2 \over \lam_1$ of the $p$--Laplacian is
bounded by a quantity depending only on the 
dimension $N$ and $p$.
\medskip
\noindent
Mathematics subject classification 35J60 35J70.
\par\noindent
keywords:  $p$--Laplacian, eigenvalue gap

\bigskip
\noindent{\bf I.  Introduction}

The linear Laplacian on a domain or a manifold can be realized as a self--adjoint 
operator, and the theory of its spectrum is a well developed subject [4,5].  For the 
$p$--Laplacian on a domain $\Omega$, 
with vanishing Dirichlet boundary conditions, 
it has been known since the work of Anane and Tsouli [1]  that a 
sequence of real eigenvalues can 
be defined by a variational procedure analogous to the min-max principle for the
linear case $p=2$,  but 
many rather basic questions about the 
spectrum remain to be addressed.  For background on 
the $p$--Laplacian, which arrives from the first variation of the functional
$${{\int_{\Omega} |\nabla u|^p} \over {\int_{\Omega} |u|^p}},\eqno(1)$$
we refer to [8, 9, 16].  
In this article $\Omega$ will be a smooth, 
bounded Euclidean domain, 
and $u$ will
vary in the Sobolev space $W_{0}^{1,p}(\Omega),$ corresponding to Dirichlet conditions on the boundary.

Several useful estimates are available in the linear case for the lowest two eigenvalues,
especially the lowest, or fundamental, eigenvalue, and some of these have been
extended to the nonlinear cases.  For instance, see [16, 17, 18].   In this article we
seek information on the ratio of the first two eigenvalues, which to our 
knowledge has not been much studied except when $p=2,$ for which case it is 
known, for example, that the
ratio is bounded universally above by the ratio attained when $\Omega$ is a ball 
[2].

We are guided by some earlier analysis of the case $p=2$ in
[11, 12, 13] as well as [2, 3], but the 
unavailability of the spectral theorem eliminates some 
essential parts of the analysis.  Therefore
we attempt to rely instead on certain
integral inequalities
as in [10].
\vskip0.3cm \noindent
Let us recall some properties of the Dirichlet $p$--Laplacian: 

\item{1.}  The $p$--Laplacian is defined for $u\in W^{1,p}(\Omega)$ by 
$\Delta_pu := \nabla \cdot \left(\left| \nabla u \right|^{p-2} \left|\nabla u
\right|\right) $. 
(Here we always equip the $p$--Laplacian with zero Dirichlet boundary
conditions.)
\item{2.}  It is then natural to define an 
eigenvalue $\lam_j$ as a value of $\lambda$ for which the eigenvalue problem
$$-\Delta_pu:= \lam|u|^{p-2}u\hbox{ in }\uom;\quad
u=0\hbox{ on }\partial\uom\eqno{\rm (2)}$$
has a nontrivial solution.  
\item{3.}  It is not, however, known whether every
such quantity is a variational eigenvalue as for the case $p=2$.
Of course it is also possible to consider the ``variational eigenvalues'' 
but it is not known in general whether these numbers coincide. 
Nevertheless, it has long been known that the first eigenvalue 
$\lambda_1$, which is isolated and simple, is the
infimum, for $u\in W^{1,p}_0(\Omega)$, of the ``Rayleigh quotient'' defined by
$(1)$.
\item{4.}  
The infimum is achieved for a multiple of a function 
$\varphi_1$ that may be chosen positive on $\Omega$. 
Henceforth we impose the normalization $\int_{\Omega} \varphi_1^p=1$. 
The pair $(\lam_1,\varphi_1)$ is referred to as the ``principal eigenpair."
\item{5.}  In 1996 Anane and Tsouli [1] gave a
characterization of the variational eigenvalues of the Dirichlet 
$p$--Laplacian 
For any positive integer $j$, let
$$\lam_j:=\inf_{{\mathcal C}\in {\mathcal C}_j}\max_{u\in {\mathcal C}}{{\int_{\Omega} |\nabla u|^p} \over {\int_{\Omega} |u|^p}},
\eqno{\rm (3)}$$
with
${\mathcal C}_j=\{{\mathcal C}\in W^{1,p}_0 (\uom)/{\mathcal C}=-{\mathcal C}; \gamma ({\mathcal C})\ge j\},$
where $\gamma$ denotes the Krasnosel'skii genus [14].   (By definition, the 
Krasnosel'skii genus of a subset ${\mathcal A}$ of
a Banach space ${\mathcal B}$ is the smallest integer $j$ for which 
there exists a nontrivial
continuous odd mapping ${\mathcal A} \rightarrow {\mathcal B}$.)
\item{6.}
Finally, it was proved by Anane and Tsouli that $\lam_2$ defined by
${\rm (3)}$ is effectively the second eigenvalue of the Dirichlet
$p$--Laplacian defined on $\uom$, in the sense that the eigenvalue
problem ${\rm (2)}$ has no other eigenvalue between
$\lam_1$ and $\lam_2$.  Hence the ratio $\lam_2 / \lam_1$ is well defined
and amenable to variational analysis. 
\vskip0.3cm 
\noindent{\bf Upper bounds on the fundamental eigenvalue ratio}
\par\nobreak
When $p=2$, bounds on $\lambda_2/\lambda_1$ can be derived from
variational estimates of the gap
$\lambda_2 - \lambda_1$.
When $p \leq 2,$ the same will be true here in certain 
situations, but 
for $p > 2$ we are led by
Lemma~3.1 in [10] to consider instead the difference
$\lam_2 - {\hat k} \lam_1$ with a suitable constant ${\hat k} \geq 1$. 
In [10] constants (depending on
$p$ and the dimension $N$)
were defined as:
$$m^p :=\max_{0\le x\le 1}((p-x)x^{p-1}+(1-x)^p), \eqno{\rm(4)}$$
and
$$\matrix{
p\geq2,\hfill&\quad N=1:   &\hat m=m =p-1 \hfill&\quad 
\hat k=k=p^{2-p}(p-1)^{p-1}\hfill\cr 
p\geq2,\hfill&\quad N\geq2:&\hat m=2^{(p-2)\over{2p}}(p-1)\hfill&\quad 
\hat k=2^{(p-2)\over2}p^{2-p}(p-1)^{p-1}\hfill\cr
1<p\leq2,\hfill&\quad N=1:   &\hat m=m \hfill&\quad \hat k = k = 1\hfill \cr
1<p\leq2,\hfill&\quad N\geq2:&\hat m=2^{(2-p)\over{2p}}m \hfill&
\quad \hat k = 1\hfill \cr }
\eqno{\rm (5)}$$
Our main result is

\proclaim Theorem 1. Let $\Omega$ be a smooth, $N$--dimensional
bounded Euclidean domain, and denote the two lowest Dirichlet
eigenvalues for the $p$--Laplacian by $\lambda_1,\lambda_2$.  

If $p \le 2$ and $N > p$, then 
$$\Gamma \,:=\, \lam_2-\lam_1 \le {{\hat m}^p N \left({p \over {N-p}}\right)^p \lam_1},
\eqno{\rm (6)}$$
or, equivalently, 
$$
{\lam_2\over\lam_1}\le \left[1+{{\hat m}^p 
\left({p \over {N-p}}\right)^p N}\right].\eqno{\rm(7)}$$
For  $p \ge 2$, 
$$\Gamma \,:=\, \lam_2-\hat k\lam_1\le {{\hat m}^p N^{-{{p} \over {2 }}} {p}^p \lam_1},
\eqno{\rm(8)}$$
or, equivalently, 
$$
{\lam_2\over\lam_1}\le \left[ \hat k+{{\hat m}^p N^{-{{p} \over {2 }}}
{p}^p}\right].\eqno{\rm(9)}$$

We prepare the proof with two estimates.  The first is a standard
(and known) uncertainty--principle inequality:

\proclaim Lemma 2.
Given a bounded domain $\Omega \in \IR^N$, with $N > p$,
for any $u \in W_0^{1,p}(\Omega)$, 
$$ \ \int_{\Omega }^{}{\left|{u  \over {\| {\bf
x}}\|_2}\right|}^{p} \leq
\left({p\over {N-p}}\right)^{p}\int_{\Omega }^{}{\left|{\nabla u }\right|}_{2}^{p},$$
where ${\bf 
x}=(x_1,\dots,x_N)$ are the Cartesian coordinates and where 
$\| {\bf x}\|_p:= (x_1^p+\cdots+x_N^p)^{1/p}$.
In particular, 
$$ \ \int_{\Omega }^{}{\left|{\varphi_1  \over {\| {\bf
x}}\|_2}\right|}^{p} \leq
\left({p\over {N-p}}\right)^{p} \lam_1, \eqno(10)$$
\medskip
\par \noindent
This follows easily from an inequality of Boggio, as shown in [10].

\proclaim Lemma 3.
Given a bounded domain $\Omega \in \IR^N$, if $p\geq 2$,
$${1 \over {\int_{\Omega}{\|{\bf x}\|_2^p \varphi_1^p}}} \leq \left(p \over N\right)^p \lam_1.$$

\noindent {\bf Proof.} 
$$\eqalign{
1=\int_\Omega\varphi^p_1&={1\over N}\int_\Omega\varphi^p_1\nabla
\cdot{\bf x}
= - {p\over N}\int_\Omega {\bf x}\cdot \varphi^{p-1}_1\nabla\varphi_1 \cr
&\le {p\over N}\left[ \int_\Omega |\nabla\varphi_1|^p\right]^{1/p}
\left[\int_\Omega \|{\bf x}\|^{p'}_2\varphi^p_1\right]^{1/p'},}$$
with $p'={p\over p-1}$ as usual, by H\"older's inequality. 
Therefore
$$1\le \left( p\over N\right)^p\lambda_1\left[
\int_\Omega \|{\bf x}\|^{p'}_2\varphi^p_1\right]^{p-1}.$$
If $p=2$, this establishes the claim.  Otherwise,  $p>2$, that is,
$p'<p$, and by H\"older's inequality, 
$$\eqalignno{
1&\le \left( p\over N\right)^p \lambda_1\left[\int_\Omega
(\|{\bf x}\|_2\varphi_1)^{p'}\varphi_1^{p(p-2)\over (p-1)}\right]^{p-1}\cr
&\le \left( p\over N\right)^p\lambda_1\left[\int_\Omega
\|{\bf x}\|^p_2\varphi^p_1\right]\left[\int_\Omega\varphi^p_1\right]^{p-2}\cr
&=\left( p\over N\right)^p\lambda_1\left[
\int_\Omega \|{\bf x}\|^p_2\varphi^p_1\right].&\square}$$

\medskip
\par \noindent
{\bf Proof of Theorem 1.}
For a unified treatment, we write $\Gamma:=\lam_2-\hat k \lam_1,$ 
noting that for $p \leq 2$,
$\hat k = 1.$
Let $\delta$ be a given real constant and set 
$\omega := \{x=(x_1,\dots,x_N) \in \Omega$, with $x_k<\delta \}$.
Assume that $\delta$ is chosen so that 
meas$(\omega)>0$ and meas$(\Omega \setminus \omega) >0$.
In other words $\delta_{\min} < \delta < \delta_{\max}$. 
Now choose $g \in {\cal C}^1_0(\Omega)$ 
such that $g(x)_{x_k=\delta}=0$.
We define
$$C:=\{\varphi_1\cdot G_{\alpha,\beta}:\alpha\in \IR, \beta\in
\IR;  |\alpha|^p+|\beta|^p=1\},$$
where
$$G_{\alpha,\beta}:=g(x)(\alpha\chi_{\omega} +
\beta\chi_{\uom \setminus \omega}).\eqno{\rm(12)}$$
It is easy to see that the set $C$ has the following properties;
\item{
(1)} $C=-C$ (change $\alpha$ to $-\alpha$ and $\beta$ to $-\beta$).
\item{
(2)} $\gamma (C) = 2$.
%
\par \noindent
  From the variational characterization ${\rm (4)}$, for any $k$, any $\delta$, and any function
$g\in {\cal C}^1_0(\Omega)$ satisfying 
$g(x)|_{x_k=\delta}=0$, it follows that
$$\lam_2\le \max_{|\alpha|^p+|\beta|^p=1}R(\alpha,\beta),\eqno{\rm(13)}$$
where
$$R(\alpha,\beta):={N(\alpha,\beta)\over D(\alpha,\beta)}$$
with
$$N(\alpha,\beta):=\int_{\uom}|\nabla (\varphi_1G_{\alpha,\beta})|^p$$
and
$$D(\alpha,\beta):=\int_{\uom}|\varphi_1G_{\alpha,\beta}|^p.$$
 From Lemma~3.1 in [10], it follows that
$$\eqalign{
N(\alpha,\beta)&\le {\hat m}^p \int_{\uom}|\nabla (G_{\alpha,\beta})
\varphi_1|^p+ \hat k \int_{\uom}{\varphi_1 |G_{\alpha,\beta}|^p
\left(-\Delta_p \varphi_1\right)}\cr
&=
{\hat m}^p \int_{\uom}|\nabla (G_{\alpha,\beta})
\varphi_1|^p+ \hat k\lam_1 \int_{\uom} |G_{\alpha,\beta}|^p\varphi^p_1.}$$
Therefore:
$$\Gamma= 
\lam_2-\hat k\lam_1 \le {\hat m}^p\max_{|\alpha|^p+|\beta|^p = 1}
{|\alpha|^p \int_\omega|\varphi_1\nabla g|^p+
|\beta|^p\int_{\uom\backslash \omega}|\varphi_1\nabla g|^p 
\over |\alpha|^p\int_\omega |\varphi_1G_{\alpha,\beta}|^p+
|\beta|^p\int_{\uom\backslash \omega}|\varphi_1G_{\alpha,\beta}|^p}\,,$$
supposing that $G_{\alpha,\beta}=0$ on $\{x_j=\delta\}$.  
The maximization in this expression is elementary:  writing $t$ for 
$|\alpha|^p$, the problem is to maximize an expression of the form
${{a t + b(1-t)} \over {c t + d(1-t)}}$ for $0 \leq t \leq 1.$  Unless this expression 
is constant and equal to ${b \over d} = {a \over c} = {{a+b} \over {c+d}}$, its 
derivative is always nonzero.  
If not constant, it is therefore
maximized when $t=0$ as $b \over d$ or else when  $t=1$ as $a \over c$.
We conclude that
$$\Gamma \le {\hat m}^p {\rm max} \left\{{{{\int_{\omega} \varphi_1^p}  \over
{\int_{\omega} \varphi^p_1|x_j-\delta|^p}}, {{\int_{\uom\backslash\omega} \varphi_1^p}  \over
{\int_{\uom\backslash\omega} \varphi^p_1|x_j-\delta|^p}}}\right\}.
\eqno{\rm(14)}$$
Observe that each of the integrals in ${\rm(14)}$  depends 
continuously on $\delta$, and that
as $\delta$ approaches the minimal value of $x_j$ in $\Omega$, 
${{\int_{\omega} \varphi_1^p}  \over
{\int_{\omega} \varphi^p_1|x_j-\delta|^p}} \rightarrow +\infty$, while 
${\int_{\uom\backslash\omega} \varphi_1^p  \over
{\int_{\uom\backslash\omega} \varphi^p_1|x_j-\delta|^p}}$
remains bounded.  The converse is the case as 
$\delta$ approaches the maximal value of $x_j$ in $\Omega$.
By continuity, there is a value of $\delta$ for which 
$${{\int_{\omega} \varphi_1^p}  \over
{\int_{\omega} \varphi^p_1|x_j-\delta|^p}} =
{{\int_{\uom\backslash\omega} \varphi_1^p}  \over
{\int_{\uom\backslash\omega} \varphi^p_1|x_j-\delta|^p}} =
{\int_{\uom} \varphi_1^p  \over
{\int_{\uom} \varphi^p_1|x_j-\delta|^p}}\,.$$
Hence

$$\int_\Omega \varphi^p_1 |x_j-\delta|^p \le {{\hat m}^p \over \Gamma}\,.$$
Let us henceforth choose the origin of the coordinate system 
so that $\delta = 0$, and then sum on $j$, obtaining
$$\int\varphi^p_1\|{{\bf x}}\|^p_p\le {{\hat m}^p \over \Gamma}\cdot N.$$
Therefore
$$\Gamma \le {{\hat m}^p N\over \int_\Omega \varphi^p_1\|{{\bf x}}\|^p_p}\,.
\eqno{\rm(15)}$$
Suppose now that $N > p$ and $p\leq 2$.  According to the Cauchy--Schwarz inequality,
$$
1
=\left(\int_{\Omega} \varphi^p_1\right)^2 = \left(\int_{\Omega} \varphi^{p/2}_1\|{\bf x}\|_p^{p/2}
\varphi_1^{p/2} \|{\bf x}\|_p^{-p/2}\right)^2\le
\int_{\Omega} \varphi^p_1\|{\bf x}\|^p_p \int_{\Omega} \varphi^p_1\|{\bf x}\|^{-p}_p.\eqno{\rm(16)}$$
We recall that  since the dimension
$N$ is finite and since $p \le 2,$
we have $ { 1 \over {\|{{\bf x}}\|_p}} \le { 1 \over {\|{{\bf x}}\|_2}}$.
\par \noindent

Hence, from ${\rm(16)}$ and ${\rm(10)}$ we derive:
$$
1\le\int_{\Omega} \varphi^p_1 \|{{\bf x}}\|_p^p \cdot \int_{\Omega}
{{\varphi_1^p}\over{\|{{\bf x}}\|_2}^p}
\le \left({p \over {N-p}}\right)^p \lambda_1 \int_{\Omega} 
\varphi^p_1\|{{\bf x}}\|_p^p\, .\eqno{\rm(17)}
$$
 Combining ${\rm(17)}$ with ${\rm(15)}$, we derive
 ${\rm(6)}$ and ${\rm(7)}$ for $p \leq 2$.

\noindent
In case $p \geq 2$, a tighter bound can be derived. 
Since for $p \ge 2,$ 
$ { 1 \over {\|{{\bf x}}\|_p} }\le N^{{p-2} \over {2 p}}
{ 1 \over {\|{{\bf x}}\|_2}}$, we deduce from Inequality ${\rm(15)}$
that
$$\Gamma \le {{\hat m}^p N^{p \over 2}\over \int_\Omega 
\varphi^p_1\|{{\bf x}}\|^p_2}\,.\eqno{\rm(18)}$$
Combining ${\rm(18)}$ with Lemma~3,
we derive  Inequalities  ${\rm(8)}$ and ${\rm(9)}$.
\hfill$\square$

It is reasonable to ask how sharp is this bound.  
Unfortunately, other than in 
the one--dimensional (or radial) case, eigenvalues of the 
$p$--Laplacian are only known 
numerically, and then essentially only the principal eigenvalue [15].  A comparison is 
possible in one-dimension, where the eigenvalues are known explicitly  
[6, 7] and
${\lambda_2 \over \lambda_1} = 2^p$.  Since the Hardy constant in one dimension is
$p \over {p-1}$, this compares to our bound of 
$p^{2-p} (p-1)^{p-1} + \left((p-1)p\right)^p$ for $p \geq 2.$  For $p=2$ it is reasonably sharp (5 rather 
than 4), but for higher values of $p$ it is less so.

As a final remark, we observe that it has been shown recently in 
[14] that for the first two eigenvalues $\Lambda_1, \Lambda_2$
of the Lindqvist $\infty$--eigenvalue
problem,  $\Lambda_k = \lim_{p \rightarrow \infty}\lambda_k(p)^{1/p}$. 
As a direct consequence of Theorem 1, we therefore have:

\proclaim Corollary 4.
$$\lim_{p \rightarrow \infty} {1 \over p} 
{\left(\lambda_2(p) \over \lambda_1(p)\right)^{1/p} }
\leq {\hat m  \over \sqrt N}.$$

$$===========================$$

\noindent
{\bf REFERENCES}
\frenchspacing
\def\bitem#1{\par\parindent=55pt\noindent\hangindent
    \parindent\hbox to\parindent {[#1]\hss}\ignorespaces}

\bitem{1}  A. Anane and N.Tsouli, 
On the second eigenvalue of the $p$--Laplacian, in:  
{\it Nonlinear Partial Differential Equations} 
(from a conference in F\`es),  
A. Benkirane and J-P Gossez, eds., 
Pitman Research Notes in Mathematics, 343.  Longman: Harlow and New York, 
1996, 1--9.


\bitem{2}  Mark S. Ashbaugh and Rafael D. Benguria, 
Proof of the Payne-P\'olya-Wein\-gerger conjecture.  
{\it Bull. Amer. Math. Soc. \bf 25} (1991) 19--29.

\bitem{3}  Mark S. Ashbaugh and Rafael D. Benguria, 
Isoperimetric inequalities for eigenvalue ratios, in 
{\it Partial Differential Equations of Elliptic Type}, Cortona, 1992, 
A. Alvino, E. Fabes, and G. Talenti, eds.  Symposia Math. 35.  
Cambridge:  Cambridge University Press, 1994, 1--36.

\bitem{4} Desmond Edmunds and Desmond Evans, {\it Spectral Theory 
and Differential Operators}, Oxford:   Clarendon, 1987.

\bitem{5}  Davies, E.B., Spectral Theory and Differential Operators, 
{\it Cambridge Studies in Advanced in Mathematics \bf 42}, Cambridge:  
Cambridge University Press, 1995.

\bitem{6}  Pavel Dr\'abek, Ranges of a--homogeneous operators and 
their perturbations, {\it \v{C}asopis 
P\v{e}st. Mat. \bf 105} (1980) 167--183.

\bitem{7} Manuel del Pino, Pavel Dr\'abek, and Raul Man\'asevich, 
The Fredholm alternative at the first eigenvalue for the one-dimensional 
$p$--Laplacian, {\it J. Diff. Eq. \bf 151} (1999) 386--419.

\bitem{8}  Pavel Dr\'abek, Pavel Krej\v{c}\'i, and 
Peter Tak\'ac, editors and authors, {\it Nonlinear Differential Equations}, 
Boca Raton, FL:  CRC Press, 1999.

\bitem{9}  Pavel Dr\'abek, A. Kufner, and F. Nicolosi, 
Nonlinear Elliptic Equations, Singular and Degenerate Case, Pilsen, 
Czech Republic, Univ. or West Bohemia, 1996

\bitem{10}  Jacqueline Fleckinger Evans M. Harrell II, and 
Fran\c{c}ois de Th\'elin,
Boundary behavior and $L^q$ estimates for solutions of equations 
containing the $p$--Laplacian, 
{\it Electronic J. Diff. Eqns. \bf 1999} (1999) No. 38, 1--19.

\bitem{11}  Evans M. Harrell II and Joachim Stubbe, On trace 
identities and universal eigenvalue estimates for some partial 
differential operators, {\it Trans. Amer. Math. Soc. \bf 349}
(1997) 1797--1809.

\bitem{12}  Evans M. Harrell II, General Bounds for the Eigenvalues 
of Schr\"odinger Operators, in {\it Maximum Principles and 
Eigenvalue Problems in Partial Differential Equations}, 
P.W. Schaefer, editor.  Essex, England:  Longman House, and New York: 
Wiley, 1988, 146--166. 

\bitem{13}  Evans M. Harrell II, Some Geometric Bounds on 
Eigenvalue Gaps, {\it Commun. in Partial Diff. Eqs. \bf 18} (1993) 179--198.

\bitem{14}  Petri Juutinen and Peter Lindqvist, On 
the higher eigenvalues for the $\infty$--eigenvalue problem, preprint 2003.

\bitem{15} Lew Lefton and Dongming Wei, Numerical approximation of 
the first eigenpair of the $p$--Laplacian using finite elements and 
the penalty method,
{\it Numer. Funct. Anal. and Optimiz. \bf 18} (1997) 389--399.

\bitem{16}  Peter Lindqvist, On a nonlinear eigenvalue problem, 
preprint 2000 (update of P. Lindqvist,
On a nonlinear eigenvalue problem).
Fall School in Analysis (Jyv\"askyl\"a, 1994), 33--54,
Report, 68, Univ. Jyv\"askyl\"a, Jyv\"askyl\"a, Finland, 1995.

\bitem{17}
Matei, Ana-Maria. First eigenvalue for the $p$--Laplace operator. 
{\it Nonlinear Anal. \bf 39} (2000), no. 8, 
Ser. A: Theory Methods, 1051--1068.

\bitem{18}
H. Takeuchi, On the first eigenvalue of the $p$--Laplacian in a 
Riemannian manifold, {\it Tokyo J. Math. \bf 21} (1998).

%
%
%
%
%
%
%

%
%

\end